%%%%%%%%%%%%%%%%%%%%%%%%%%%%%%%%%%%%%%%%%%%%%%%%%%%%%%
%%%%%%%%%%% TOP OF TEX FILE %%%%%%%%%%%%%%%%%%%%%%%%%%
%%%%%%%%%%%%%%%%%%%%%%%%%%%%%%%%%%%%%%%%%%%%%%%%%%%%%%
%%%%%%%%%%% Definitions     %%%%%%%%%%%%%%%%%%%%%%%%%%
\magnification=\magstephalf
\parindent=0mm
\parskip=0mm
\hsize=5.75 true in
\baselineskip=15pt 
\bigskipamount=24pt plus 6 pt minus 6 pt
\medskipamount=18pt plus 3 pt minus 3 pt
\smallskipamount=12pt plus 2 pt minus 2 pt
\thinmuskip=5mu
\medmuskip=6mu plus 3mu minus 2mu
\thickmuskip=6mu plus 3mu minus 2mu
\overfullrule=0pt
\font\xivrm=cmr17
%mathemath. Zeichen:

\def\R{{I\mkern-5mu R}}
\def\N{{I\mkern-5mu N}}

\def\Z{{Z\mkern-8mu Z}}

\def\Any{\,\cdot\,}

\def\norm#1{\mathopen\Vert #1 \mathclose\Vert}
\def\Modul#1{\left\vert #1 \right\vert}
\def\modul#1{\mathopen\vert #1 \mathclose\vert}

\def\sgn{\mathop{\rm sgn }}

\def\loc{\sb{\rm loc}}

\def\square{\mathord{\vbox{\hrule\hbox{\vrule\hskip 9pt\vrule height
9pt}\hrule}}}
\frenchspacing
%%%%%%%%%%%%%%%%%%%%%%%%%%%%%%%%%%%%%%%%%%%%%%%%%%%%%%%%%%%%%%%%%
%%%%%%%%%%%  Beginning of main text  %%%%%%%%%%%%%%%%%%%%%%%%%%%%
%%%%%%%%%%%%%%%%%%%%%%%%%%%%%%%%%%%%%%%%%%%%%%%%%%%%%%%%%%%%%%%%%
{\xivrm
\centerline{Spherically symmetric Dirac operators with variable}
\centerline{mass and potentials infinite at infinity.}
}

\bigskip
\line{\hfill{\it Karl Michael Schmidt$^1$ and Osanobu Yamada$^2$}\hfill}
\medskip
\line{$^1$Mathematisches Institut der Universit\"at, Theresienstr. 39,
      D-80333 M\"unchen, Germany}
\line{$^2$Department of Mathematics, Ritsumeikan University, Kusatsu,
      Shiga 525-77, Japan}

\bigskip
1991 Mathematics Subject Classification: 81Q10, 35Q40, 34L40

\bigskip
{\bf Abstract.}\quad
{\it We study the spectrum of spherically symmetric Dirac operators
in three-dimensional space with potentials tending to infinity
at infinity under weak regularity assumptions.
We prove that purely absolutely continuous spectrum covers the
whole real line if the potential dominates the mass, or scalar
potential, term.
In the situation where the potential and the scalar potential
are identical, the positive part of the spectrum is purely
discrete; we show that the negative half-line is filled with
purely absolutely continuous spectrum in this case.}

\bigskip
{\bf \S1.\ \ Introduction.}

\smallskip
In a recent paper [20] the spectral properties of the
three-dimensional Dirac operator
$$
  H = \alpha\cdot p + m(x)\,\beta + q(x)\,I_4
  \qquad (x\in\R^3)
$$
(where $p = -i\nabla$,
$i^2 = -1$, $I_n$ is the $n\times n$ unit matrix,
and $\alpha_0 = \beta, \alpha_1, \alpha_2, \alpha_3$
are Hermitian $4\times 4$ matrices
satisfying the anti-commutation relations
$$
  \alpha_j\,\alpha_k + \alpha_k\,\alpha_j = 2 \delta_{jk}
  \qquad (j,k \in \{ 0, 1, 2, 3 \}))
$$
were studied under the condition that
the real-valued coefficient function $m$
tends to $\infty$ (or $-\infty$) as $\modul{x}\rightarrow\infty$.
For constant $m$, $H$ is the Hamilton operator describing a relativistic
quantum mechanical particle of mass $m$ moving in an external force field
of (real-valued) potential $q$.
As a non-constant function, $m$ can also take the role of a so-called
scalar
potential, which has been discussed in the physical literature as a model
of
quark confinement (cf. the references in [20], Thaller [14] p. 305).

In [20] it is shown that if $m$ dominates $q$ and tends to infinity
as $\modul{x} \rightarrow\infty$, then the spectrum of the Dirac operator
$H$
is purely discrete;
if $m$ coincides identically with $q$,
and tends to infinity (in addition to certain regularity requirements),
then the positive part of the spectrum of $H$ is purely discrete.
Furthermore, if $m\equiv q$ is of at most quadratic growth, the negative
half-line is filled with purely continuous spectrum of $H$.

It seems to be a rather more delicate question to determine the quality
(absolute continuity or otherwise) of the continuous part of the spectrum
of
$H$ in the case $m \equiv q$, or even the overall structure of the spectrum
in the situation where $q$ dominates $m$ and tends to infinity.
It is the purpose of the present paper to address this question under the
additional assumption that $m$ and $q$ are spherically symmetric functions;
then the operator $H$ is spherically symmetric in the sense that rotations
in space lead to unitarily equivalent operators.

By the well-known procedure of separation in spherical polar coordinates
(cf. Weidmann [18] Appendix to Section 1),
$H$ is then unitarily equivalent to the direct sum of
the countable family of one-dimensional Dirac operators on the half-line
$r > 0$,
$$
  h_k = \sigma_2\, p + {k \over r}\,\sigma_1 + m(r)\,\sigma_3 + q(r)\,I_2
  \qquad (k \in \Z\setminus\{0\}),
$$
where ${\displaystyle p = -i\,{d\over dr}}$,
and $\sigma_1 = \pmatrix{0 & 1 \cr 1 & 0 \cr}$,
$\sigma_2 = \pmatrix{0 & -i \cr i & 0 \cr}$,
$\sigma_3 = \pmatrix{1 & 0 \cr 0 & -1 \cr}$
are the Pauli matrices.
If $m$ and $q$ are regular, i.e. locally integrable on $(0, \infty)$ and
integrable at $0$, this Dirac system is in the limit point case at $0$ for
all $k \in \Z\setminus\{0\}$ by virtue of the singular angular momentum
term $k\,\sigma_1/r$ ([11] Lemma 1);
moreover, it is in the limit point case at $\infty$ as well
(Weidmann [18] Corollary to Theorem 6.8), and consequently the minimal
operator
associated with the formal expression $h_k$ is essentially self-adjoint
by the Weyl theory (Weidmann [18] Theorem 5.8).
We denote the unique self-adjoint realization again by $h_k$;
$H := \bigoplus_{k\in\Z\setminus\{0\}} h_k$ is a self-adjoint realization
of the three-dimensional Dirac operator.

The study of the spectrum of $H$ can then essentially
be reduced to that of the individual one-dimensional operators $h_k$.
The spectral properties of half-line Dirac operators with potentials which
do not approach a finite limit at infinity, have been studied previously 
in the literature in different special situations.
For example, Hinton and Shaw [8], extending the work of Roos and Sangren
[9], give
conditions for a Dirac system with dominant $m$ to have purely discrete
spectrum.
On the other hand,  Evans and Harris [4] derive lower bounds on the
absolute value of eigenvalues in a situation with dominant $q$.
In the case of the electron Dirac operator ($m = \hbox{const}$), Erd\'elyi
[3] has established that the spectrum of $h_k$ is purely absolutely
continuous and covers the whole real line provided $q$ is locally
absolutely
continuous,
$$
  \lim_{r\rightarrow\infty} q(r) = \infty,
\qquad
  \hbox{and}
\qquad
  \int^\infty {\modul{q'} \over q^2} < \infty,
$$
thus refining a result of Titchmarsh [15] that had been anticipated,
in a formal way, by Rose and Newton [10].
Recently a different proof of this assertion (for $m \equiv 1$), based on
the Gilbert-Pearson method, was given by [12] under the slightly weaker
hypotheses
$q\in BV\loc(0,\infty)$,
$$
  \lim_{r\rightarrow\infty} q(r) = \infty,
\qquad
  \hbox{and}
\qquad
  {1 \over q} \in BV(\Any,\infty)\,.
$$
(Here $BV(I)$ denotes the space of functions of bounded variation on the
interval $I \subset \R$, and $BV\loc(I)$ is the corresponding local space.
Generally, if $X(I)$ is a space of functions on the real interval $I$, we
here and below use the notation
$$
  X(\Any,\infty) := \{f\mathop{|} \hbox{there is } a\subset\R \hbox{ such
that }
    f\in X([a,\infty))\}.)
$$

In this paper we consider, following the methods of [12], the case of
non-constant $m$; this also requires that the angular momentum term be
handled
in a different way than the perturbative treatment indicated in [12],
(cf. Remark 6 below).
We prove that (under certain weak regularity assumptions) the negative
part of the spectrum of $H$ is purely absolutely continuous if
$$
  m(r) \equiv q(r) \rightarrow \infty \qquad (r \rightarrow\infty),
$$
and that the whole real line is filled with purely absolutely continuous
spectrum of $H$ if
$$
  m(r) << q(r) \rightarrow \infty \qquad (r \rightarrow\infty).
$$

\bigskip 
{\bf \S2.\ \ Results.}
\nobreak
\smallskip
We prove the following theorem in Section 3.

\medskip
{\bf Theorem 1.}\qquad
{\it
Let $q \in L^1\loc([0, \infty))$, $m\in AC\loc([0, \infty))$, and assume
that
$$\eqalign{
  \hbox{\rm (A1)}&\qquad \lim_{r\rightarrow\infty} q(r) = \infty,
\cr
  \hbox{\rm (A2)}&\qquad \liminf_{r\rightarrow\infty}\, \modul{m(r)} > 0,
  \qquad \limsup_{r\rightarrow\infty} \Modul{{m(r) \over q(r)}} < 1,
\cr
  \hbox{\rm (A3)}&\qquad {m \over q - \lambda} \in BV(\Any,\infty) \qquad
(\forall \, \lambda\in\R),
\cr
  \hbox{\rm (A4)}&\qquad {m' \over r m q} \in L^1(\Any,\infty).
\cr} $$
Then $\sigma_{ac} (H) = \R$, and $\sigma_s (H) = \emptyset$.
}

\medskip
{\it Remark 1.}\qquad
By unitary equivalence, the same result holds true if, instead of (A1),
we assume $\lim_{r\rightarrow\infty} q(r) = - \infty$. In view of the
condition (A3) we may rewrite the latter condition in (A2) as follows ;
$$ \lim_{r \rightarrow \infty} \Modul{{m(r) \over q(r)}} < 1.$$
\smallskip
{\it Remark 2.}\qquad
If $m$ is constant, (A3) is true for all $\lambda \in \R$ if it is
true for one real value of $\lambda$ (cf. Remark 8 in the Appendix).
This does not hold in the general case; for example, choosing
$m(r) := 2 + \sin r$, $q(r) = r^{1/4} m(r)$ ($r\in [0,\infty)$),
we find that (A1), (A2), (A4) and (A3) for $\lambda = 0$ are satisfied,
but (A3) for $\lambda \neq 0$ is not.

However, Proposition 5 (in the Appendix) shows that it is generally
sufficient
to assume (A3) for {\it two} distinct values of $\lambda$.
Alternatively, if we assume, in addition to (A1),
that $q \in AC\loc(\Any,\infty)$ and
$m q' q^{-3} \in L^1(\Any,\infty)$,
then
$m / q \in BV(\Any, \infty)$ implies
$m / (q - \lambda) \in BV(\Any, \infty)$ for all
$\lambda \in \R$, by a straightforward application of Proposition 3 in the
Appendix to
${\displaystyle
  {m \over q - \lambda} = {q \over q - \lambda} \cdot {m \over q}.
}$

\medskip
Theorem 1 (cf. also Theorem 2 below),
as well as the results of [12], Weidmann [16] \S 8, and Weidmann [17],
strongly suggest that a condition of bounded variation of the
coefficient functions is a natural setting for statements on absolute
continuity of the spectrum of one-dimensional Hamiltonians.
Nevertheless, assuming higher regularity of the coefficients (as will
usually be given in concrete applications) one can reformulate the
conditions of bounded variation in terms of integrability of derivatives,
yielding hypotheses which may be more convenient to verify than (A3)
itself.
Thus we note

\medskip
{\bf Corollary 1.}\qquad
{\it
Let $q, m \in AC\loc([0,\infty))$ and assume {\rm (A1)}, {\rm (A2)}, and
$$
  {m' \over q}, \quad {mq' \over q^2} \in L^1(\Any, \infty).
$$
Then $\sigma_{ac}(H) = \R$, $\sigma_s(H) = \emptyset$.
}

\medskip
{\it Proof.}\qquad
The assumption gives
${\displaystyle
  \left({m \over q - \lambda}\right)'
  = {m'q - q'm - \lambda m' \over (q - \lambda)^2}
  \in L^1(\Any, \infty)
}$
for each $\lambda \in \R$,
which implies (A3).
Furthermore, by (A2) we have $r \modul{m(r)} \rightarrow \infty$
($r \rightarrow \infty$), and thus
$$
  {m' \over r m q} = {1 \over r m} \cdot {m' \over q} \in L^1(\Any,
\infty).
\eqno\square
$$

\bigskip
As a borderline case between the two situations in which either $q$
dominates $m$ (giving rise to purely absolutely continuous spectrum as
shown above), or $m$ dominates $q$ (with purely discrete spectrum as a
consequence of [20]), we then study the case in which $m$ and $q$ coincide
identically.
In Section 4 we prove the following result.

\medskip
{\bf Theorem 2.}\qquad
{\it
Let $m \equiv q \in AC\loc([0, \infty))$, and assume
$$\eqalign{
  \hbox{\rm (B1)}&\qquad \lim_{r\rightarrow\infty} q(r) = \infty,
\cr
  \hbox{\rm (B2)}&\qquad {q' \over q^{3/2}} \in BV(\Any, \infty)
                   \cap L^2(\Any, \infty).
\cr}$$
Then $\sigma_{ac}(H) = (-\infty, 0]$,
$\sigma_s(H) \cap (-\infty, 0) = \emptyset$.
}

\medskip
For the convenience of the reader, we include a corollary in which
assumption
(B2), which involves bounded variation, is replaced by a condition of
integrability of derivatives.

\medskip
{\bf Corollary 2.}\qquad
{\it
Let $m \equiv q \in C^2([0,\infty))$, and assume {\rm (B1)} and
$$
  \hbox{\rm (B2)'}\qquad {q''\over q^{3/2}},
  \quad {(q')^2 \over q^{5/2}} \in L^1(\Any,\infty).
$$
Then $\sigma_{ac}(H) = (-\infty, 0]$, $\sigma_s(H) \cap (-\infty, 0) =
\emptyset$.
}

\medskip
{\it Proof.}\qquad
${\displaystyle
  \left({q' \over q^{3/2}}\right)'
= {q'' \over q^{3/2}} - {3 (q')^2 \over 2 q^{5/2}} \in L^1(\Any,\infty)
}$
and

${\displaystyle
  \left({q' \over q^{3/2}}\right)^2
= {1 \over q^{1/2}}\cdot {(q')^2 \over q^{5/2}}
\in L^1(\Any,\infty)
}$
imply (B2).
\hfill $\square$

\medskip
{\it Remark 3.}\qquad
Under the hypotheses of Corollary 2, every solution $u$ of
$h_k\,u = \lambda\,u$
is twice continuously differentiable, and satisfies the ordinary
differential
equation system
$$
  u_1' + {k \over r}\,u_1 = \lambda\,u_2,
\qquad
  -u_1'' + 2 \lambda q\,u_1 + {k(k + 1) \over r^2}\,u_1 = \lambda^2\,u_1.
$$
Along the lines of Dunford and Schwartz [2] Theorem XIII.6.20,
one can then prove the existence of a fundamental
system $u_+, u_-$ with asymptotic behaviour
$$
  u_\pm (r) = \pmatrix{
   q(r)^{-1/4} \left(e^{\pm i \int_1^r \sqrt{\lambda^2 - 2 \lambda q}} +
o(1)
               \right)
\cr
   \mp i \sqrt{{2 \over -\lambda}}\,q(r)^{1/4}
               \left(e^{\pm i \int_1^r \sqrt{\lambda^2 - 2 \lambda q}} +
o(1)
               \right)
\cr}
$$
as $r \rightarrow \infty$, which implies the non-existence of subordinate
solutions for negative $\lambda$.
In order to achieve weaker regularity requirements, our proof of Theorem 2
follows a different approach, remaining fully in the context of Dirac
systems
without recourse to associated second order differential equations.

\bigskip 
{\bf \S3.\ \ The\ \ Case\ \ $m << q$.}
\nobreak
\smallskip
In this Section we prove Theorem 1.
We proceed as follows.
As indicated in the Introduction, the spherically symmetric Dirac operator
$H$
is unitarily equivalent to the countable direct sum of half-line operators
with angular momentum term:
$$
  H \cong \bigoplus_{k\in \Z\setminus\{ 0\}} h_k;
$$
therefore we have
$$
  \sigma_{ac}(H) = \overline{\bigcup_{k\in\Z\setminus\{0\}}
\sigma_{ac}(h_k)},
\qquad
  \sigma_s(H) = \overline{\bigcup_{k\in\Z\setminus\{0\}} \sigma_s(h_k)}
$$
(cf. [11]), and hence it is sufficient to prove that for each nonzero
integer $k$, $h_k$ has purely absolutely continuous spectrum throughout
the real line.

By the following proposition, this can be reduced to showing that all
solutions of the eigenvalue equation for $h_k$, for all real values of
the spectral parameter, are bounded at infinity.

\medskip
{\bf Proposition 1.}\qquad
{\it
Let $a \ge -\infty$, $l, m, q \in L^1\loc((a, \infty))$ be real-valued
functions, and $I\in\R$ an open interval.
Let $h$ be a self-adjoint realization of the Dirac system
$$
  \sigma_2\,p + m\,\sigma_3 + l\,\sigma_1 + q\,I_2
$$
on $(a, \infty)$.
If, for every $\lambda \in I$, every solution $u$ of
$$
  (\sigma_2\,p + m\,\sigma_3 + l\,\sigma_1 + q\,I_2)\,u = \lambda\,u
$$
is an element of $L^\infty(\Any, \infty)$,
then we have $I \subset \sigma_{ac}(h)$, $\sigma_s(h) \cap I = \emptyset$.
}

\medskip
{\it Remark 4.}\qquad
It is possible to give a direct proof for this special case of the
Gilbert-Pearson method for Dirac systems (Behncke [1] Lemma 1) along the
lines of Simon [13], cf. the Appendix of [12], bypassing the complications
of the general theory.
Note that only a condition for each separate value of $\lambda$ is
imposed in the above Proposition;
however, in our situation, one could also apply Weidmann's theorem
[19] instead of Proposition 1, after showing that the bounds on the
solutions are locally uniform in $\lambda$.

\medskip
The heart of the matter is contained in Proposition 2, in which a
sufficient
condition for the boundedness of all solutions of a general Dirac
system
$$
  (\sigma_2\,p + M\,\sigma_3 + L\,\sigma_1 + Q\,I_2)\,u = 0,
\eqno (*)
$$
with $\lim_{r\rightarrow\infty} Q(r) = \infty$, is given.

After proving Proposition 2, we finish the Proof of Theorem 1 by showing
that
$Q := q - \lambda$, $M := m$ and $L := k/r$ satisfy the hypotheses of
Proposition 2 for all real $\lambda$ and all non-zero integers $k$.

\medskip
{\bf Proposition 2.}\qquad
{\it
Let $a \ge -\infty$, $Q, M, L \in L^1\loc(a, \infty)$ be real-valued
functions such that
$$\eqalign{
  \hbox{\rm (C1)}&\qquad \lim_{r\rightarrow\infty} Q(r) = \infty,
\cr
  \hbox{\rm (C2)}&\qquad \limsup_{r\rightarrow\infty} {W(r) \over Q(r)}
      < 1,
\cr
  \hbox{\rm (C3)}&\qquad {W \over Q - W},\,
              {M \over Q - W},\,
              {L \over Q - W} \in BV(\Any, \infty),
\cr}$$
where $W := \sqrt{M^2 + L^2}$.
Then every solution of $(*)$ is bounded at infinity.

If $M \equiv 0$ {\rm [}$L \equiv 0${\rm ]},
condition {\rm (C3)} can be replaced by
$$
  \hbox{\rm (C3)'} \qquad {L \over Q - L} \in BV(\Any,\infty)
\qquad
  \left[
  {M \over Q - M} \in BV(\Any,\infty)
  \right].
$$
}

\medskip
{\bf Corollary 3.}\qquad
{\it
Under the hypotheses of Proposition 2, let $r_0 > a$.
Then there is a constant $C > 0$ such that
$$
  {\modul{y(r_0)^2} \over C} \le \modul{y(r)}^2 \le C\,\modul{y(r_0)}^2
\qquad (r \ge r_0)
$$
holds for every solution $y$ of $(*)$.
}

\medskip
{\it Proof.}\qquad
By Proposition 2 every fundamental system of $(*)$ is bounded on
$[r_0,\infty)$.
Consequently, there is $C > 0$ such that
$$
  \modul{y(r)}^2 \le C\,\modul{y(r_0)}^2
  \qquad (r \ge r_0)
$$
holds for all solutions $y$ of $(*)$
Let $z$ be the solution of $(*)$ with initial value
$z_1(r_0) = -y_2(r_0)$, $z_2(r_0) = y_1(r_0)$;
as the Wronskian is constant, we have
$$
  \modul{y(r_0)}^4
  = \left|\matrix{y_1(r_0) & -y_2(r_0) \cr y_2(r_0) & y_1(r_0)
\cr}\right|^2
  = \left|\matrix{y_1 & z_1 \cr y_2 & z_2 \cr}\right|^2
  \le \modul{y}^2\modul{z}^2
  \le C\,\modul{y}^2\modul{z(r_0)}^2,
$$
and the assertion follows.
\hfill $\square$

\medskip
{\it Remark 5.}\qquad
If we specialize $L \equiv 0$, $M\equiv 1$, Proposition 2 states that
for all real $\lambda$, all solutions of the differential equation
$$
  (\sigma_2\,p + \sigma_3 + q\,I_2)\,u = \lambda u
$$
are bounded at infinity provided $\lim_{r\rightarrow\infty} q(r) = \infty$
and $1/(q - \lambda)\in BV(\Any,\infty)$ (which is equivalent to $1/q \in
BV(\Any,\infty)$, cf. Remark 8 in the Appendix).
These assumptions are equivalent to the hypotheses of Theorem 1 of [12],
since a bounded function of locally bounded variation has bounded total
variation if and only if it has bounded positive variation.

\medskip
{\it Remark 6.}\qquad
Conversely, seeing that, in contrast to the rather large mass term
$\sigma_3\,m$, the angular momentum term $(k/r)\sigma_1$ in the eigenvalue
equation for $h_k$ is a smooth function on $(0, \infty)$
decaying at infinity,
one may at first be tempted to treat it as
a perturbation of the equation without angular momentum term, which can be
handled by a simpler version of Proposition 2 (this line of attack was
sketched in [12] Remark 2, for the case $m = 1$).
However, the unitary transformation used to turn the angular momentum
term into an integrable perturbation ([11] Lemma 3), involves the
derivatives
of both the angular momentum and the mass (or scalar potential) term,
leading to the stronger requirement
$$
  \hbox{\rm (A4)'} \qquad {m' \over r m^2} \in L^1(\Any,\infty)
$$
instead of (A4) in Theorem 1.
(For example, if $q(r) > c\, r^\varepsilon$ for some $c, \varepsilon > 0$,
and $m$ is periodic, then (A4) is satisfied, but (A4)' is not.)
Therefore we prefer this strong version of Proposition 2, in which $L$ and
$M$ enter in a perfectly symmetric way.
On quite different grounds, Proposition 2 will also prove useful in Section
4.

\medskip
The basic idea for the proof of Proposition 2 is essentially that developed
in
the proof of Theorem 1 of [12].
Instead of the pointwise norm
$\modul{u} := \sqrt{\modul{u_1}^2 + \modul{u_2}^2}$
of the solution itself, we study the behaviour of an associated quantity
$R(r)$, which may be interpreted as the major radius of the elliptical
orbit in the $(u_1, u_2)$ plane on which the (real-valued) solution would
be running if the coefficients were held constant at their value at $r$.
In the case $L \equiv 0$ considered in [12] the major axes of the
ellipse coincide with the $u_1$-, $u_2$-axes;
in the general situation captured by Proposition 2, however, the ellipse is
oblique, which renders the function $R$, as expressed in terms of the
solution $u$, considerably more complicated. ---
Note that our $R$ corresponds to the $R^2$ of [12].

\medskip
{\it Proof of Proposition 2.}\qquad
Let $u \in AC\loc([a, \infty))$ be a solution of $(*)$ with real-valued
components.
By hypothesis, there is $r_0 > 0$ such that
  $\quad Q(r) - W(r) > 0 \quad (r \ge r_0)$, and
$$
  {W \over Q - W},\quad
  {M \over Q - W},\quad
  {L \over Q - W} \in BV([r_0, \infty)).
$$
On $[r_0, \infty)$ we consider the function
$$
  R := {1 \over Q - W}\,
       ((u_1^2 + u_2^2)\, Q + (u_1^2 - u_2^2)\, M + 2\,u_1\,u_2\, L).
$$
Note that
$$
  {Q \over Q - W}
  = 1 + {W \over Q - W}
  \in BV([r_0, \infty)).
$$
Thus we have $R \in BV\loc([r_0, \infty))$ by Proposition 3 in the
Appendix.

For $t_2 \ge t_1 \ge r_0$, the formula for integration by parts for
Stieltjes
integrals (Fichtenholz [5] \S 577) and the local absolute continuity of $u$
yield:
$$\eqalign{
  R(t_2) - R(t_1)
  = \int_{t_1}^{t_2} &\left( (u_1^2 + u_2^2)'\, Q + (u_1^2 - u_2^2)'\, M
                            + (2\, u_1\, u_2)'\, L \right)
      \, {1 \over Q - W}
\cr
  &+ \int_{t_1}^{t_2} (u_1^2 + u_2^2)\,
         d\left({W \over Q - W}\right)
  + \int_{t_1}^{t_2} (u_1^2 - u_2^2)\,
         d\left({M \over Q - W}\right)
\cr
  &+ \int_{t_1}^{t_2} 2\, u_1\, u_2\,
         d\left({L \over Q - W}\right).
\cr}$$

As $u$ is a solution of $(*)$ and hence
$$
  u_1' = - L\, u_1 - (Q - M) \, u_2,
\qquad
  u_2' = (Q + M)\, u_1 + L\, u_2,
$$
a straightforward calculation shows that the integrand of the first
integral vanishes identically.
Estimating
$\modul{u_1^2 - u_2^2}, 2\, \modul{u_1\, u_2} \le \modul{u}^2$, and
$$
  R = u_1^2 + u_2^2
    + {1 \over Q - W} \left(
       \sqrt{W + M}\, u_1 +
       (\sgn L)\, \sqrt{W - M}\, u_2 \right)^2
    \ge \modul{u}^2,
$$
and noting that
${\displaystyle\quad
  \mathop{\rm Var}_{[t_1, t_2]} f = \int_{t_1}^{t_2} \modul{df}
}$\quad
holds for every function $f$ of locally bounded variation, we find
$$
  R(t_2) - R(t_1)
  \le \left(
\mathop{\rm Var}_{[t_1, t_2]} {W \over Q - W} +
\mathop{\rm Var}_{[t_1, t_2]} {M \over Q - W} +
\mathop{\rm Var}_{[t_1, t_2]} {L \over Q - W} \right)\,
\sup_{[t_1, t_2]} R.
$$

By Lemma 2 of [12], it follows that $R$, and consequently $u$, is
bounded in $[r_0, \infty)$.
This concludes the proof in the general case.

In the case $M \equiv 0$, we consider the function
$$
  R := u_1^2 + u_2^2 + {L \over Q - L}\,(u_1 + u_2)^2.
$$
By virtue of (C2) and (C3)' there is some $\varepsilon \in (0, 1)$
and $r_0 > 0$ such that $\modul{L(r)} \le \varepsilon\, Q(r)$ $(r \ge
r_0)$,
and $L / (Q - L) \in BV([r_0, \infty))$.
Therefore $R \in BV\loc([r_0, \infty))$, and it is not difficult to check
that
$$
  {2 L(r) \over Q(r) - L(r)} \ge {-2\varepsilon \over 1 + \varepsilon}
  > - 1,
$$
and hence ${\displaystyle R(r) \ge {1 - \varepsilon \over 1 + \varepsilon}
\modul{u(r)}^2}$
$(r \ge r_0)$.
By essentially the same calculation as above we infer that
for $t_2 \ge t_1 \ge r_0$,
$$
  R(t_2) - R(t_1) \le 2\,{1 + \varepsilon \over 1 - \varepsilon}\,
    \sup_{[t_1, t_2]} R \cdot \mathop{\rm Var}_{[t_1, t_2]} {L \over Q -
L},
$$
and the assertion again follows by Lemma 2 of [12].

For $L \equiv 0$, take
$$
  R := u_1^2 + u_2^2 + {2 M \over Q - M}\,u_1^2,
$$
and proceed as above.
\hfill $\square$

\medskip
Now we conclude the proof of Theorem 1, showing that under hypotheses
(A1) -- (A4), $Q := q - \lambda$, $M := m$ and $L(r) := k/r$ $(r > a := 0)$
satisfy (C1) -- (C3), for every $\lambda\in\R$ and $k \in
\Z\setminus\{0\}$.

(C1) and (C2) are obvious from (A1) and (A2), respectively.
As in Proposition 2, abbreviate $W := \sqrt{M^2 + L^2}$.
We first prove $W / (Q - W) \in BV(\Any,\infty)$,
which is equivalent to $W / Q \in BV(\Any, \infty)$ by
Proposition 4 in the Appendix.

By Proposition 3, $W / Q = (W/M)\cdot (M/Q) \in BV(\Any,\infty)$,
since $M/Q \in BV(\Any,\infty)$ by (A3),
$W/M \in L^\infty(\Any,\infty) \cap AC\loc(\Any,\infty)$, and
$$
  {M \over Q}\cdot \left({W \over M}\right)'
  = {L\,(ML' - LM') \over Q M W}
  = - {k^2\, q \over r W Q} \left({m' \over q m r} + {1 \over q r^2}\right)
  \in L^1(\Any,\infty).
$$

Similarly, we have $M /(Q - W) = (M/W)\cdot(W / (Q - W)) \in
BV(\Any,\infty)$,
since $M/W \in L^\infty(\Any,\infty) \cap AC\loc(\Any,\infty)$, and
$$
  {W \over Q - W} \cdot \left({M \over W}\right)'
  = {M Q \over W\,(Q - W)} \cdot {L\,(M'L - ML') \over Q M W}
  \in L^1(\Any,\infty).
$$

Finally, $L / (Q - W) = (L/M)\cdot (M/(Q - W)) \in BV(\Any,\infty)$,
since $L/M \in L^\infty(\Any,\infty) \cap AC\loc(\Any,\infty)$, and
$$
  {M \over Q - W}\cdot \left(L \over M\right)'
  = - {k \over r^2\,(Q - W)}
    - {k\,q \over Q - W} \cdot {m' \over q m r}
  \in L^1(\Any,\infty).
\eqno \square
$$

\bigskip
{\bf \S 4.\ \ The\ \ Case\ \ $m \equiv q$.}
\nobreak

\smallskip
In this section we prove Theorem 2.
Since it is shown in [20] that $\sigma(H) \cap (0, \infty)$ is
purely discrete, and since the spherically symmetric operator $H$ admits
separation in spherical polar coordinates, it is sufficient to prove that,
for each non-zero integer $k$, $(-\infty, 0) \subset \sigma_{ac}(h_k)$,
and $\sigma_s(h_k) \cap (-\infty, 0) = \emptyset$.

We shall make use of the Gilbert-Pearson theory
(Gilbert and Pearson [6], Gilbert [7], Behncke [1]),
showing that for negative $\lambda$, the differential equation
$$
  (\sigma_2\,p + q(r)\,\sigma_3 + q(r)\,I_2 + {k \over r}\,\sigma_1)\,u
  = \lambda\, u
\eqno (**)
$$
does not possess a subordinate solution at $\infty$, i.e.
that any two non-trivial solutions $v$ and $w$ of $(**)$, $\lambda < 0$,
satisfy
$$
  \liminf_{r\rightarrow\infty}
     {\int_{r_0}^r \modul{v}^2 \over \int_{r_0}^r \modul{w}^2} > 0
$$
for some $r_0 > 0$.
We have to go back to this more general definition, because unlike the
situation of Section 3, the solutions will be unbounded, as suggested
by the asymptotics given in Remark 3.

By a transformation which
takes into account this expected growth, or decay, behaviour of the
solutions,
we obtain a differential equation which is again of Dirac type,
with coefficients which (with the help of Lemma 1 below) can be shown
to satisfy the hypotheses of Proposition 2.
As in Section 2, we can therefore conclude that all solutions of the
transformed equation are globally bounded at infinity;
indeed Corollary 3 shows that all solutions of the transformed equation are
of the same size.
This does not immediately imply that the original equation has no
subordinate
solutions, since the transformation is unbounded;
yet considering the oscillatory behaviour of the solutions,
it turns out that the growing and decaying components are sufficiently
well distributed between solutions to prevent the
existence of a subordinate solution.

\medskip
{\bf Lemma 1.}\qquad
{\it
Under the hypotheses of Theorem 2, let $\lambda \in \R$ and set
$\gamma := 2\,q - \lambda$.
Then $\gamma \in AC\loc([0, \infty))$, and
$$
  {\gamma^{\,\prime} \over \gamma^{3/2}} \in BV(\Any,\infty),
\qquad
  {\gamma^{\,\prime} \over r\,\gamma^{3/2}} \in L^1(\Any,\infty),
\qquad
  \lim_{r\rightarrow\infty} {\gamma^{\,\prime}(r) \over \gamma(r)^{3/2}} =
0.
$$
}

\medskip
{\it Proof.}\qquad
There is $r_0 > 0$ such that $\gamma(r) \ge q(r) \ge 1$ $(r \ge r_0)$,
and $g := 2 q^{\,\prime} q^{-3/2} \in BV([r_0, \infty)) \cap L^2([r_0,
\infty))$.

Applying Proposition 3 in the Appendix with $f := (q / \gamma)^{3/2}$
(noting that $\modul{f} \le 1$,
$$
  \modul{f^{\prime} g} = \Modul{{3\lambda q^{\prime 2} \over q
\gamma^{5/2}}}
  \le 3 \modul{\lambda} {q^{\prime 2} \over q^{7/2}}
  \in L^1([r_0, \infty))),
$$
we obtain $\gamma^{\,\prime} \gamma^{-3/2} = f g \in BV([r_0, \infty))$.

Furthermore,
$\modul{\gamma^{\,\prime} \gamma^{-3/2}} \le \modul{g}$,
and thus the second assertion follows by the Schwarz inequality.
Finally, the last assertion follows from the fact that $g$, as a
square integrable function of bounded variation on $[r_0, \infty)$, must
converge to $0$ at infinity.
\hfill $\square$

\medskip
{\it Proof of Theorem 2.}\qquad
Let $\lambda < 0$, $k \in \Z \setminus \{0\}$, and $u$, $y$ nontrivial
solutions of $(**)$ with real-valued components.
On a right half-axis on which $\gamma : = 2\, q - \lambda$ is positive,
we consider the functions
$$
  v := \pmatrix{ (\gamma / \Lambda)^{1/4}\, u_1 \cr
                 (\Lambda / \gamma)^{1/4}\, u_2 \cr},
\qquad
  w := \pmatrix{ (\gamma / \Lambda)^{1/4}\, y_1 \cr
                 (\Lambda / \gamma)^{1/4}\, y_2 \cr}
\qquad
  (\hbox{with } \Lambda := \modul\lambda);
$$
then $v$ and $w$ are solutions of the differential equation of Dirac type
$$
  (\sigma_2\,p + L\,\sigma_1 + Q\,I_2)\, v = 0,
$$
where ${\displaystyle L = {k \over r} - {\gamma^{\,\prime} \over 4
\gamma}}$,
$Q = \sqrt{\Lambda\gamma}$.

We observe that $\lim_{r\rightarrow\infty} Q(r) = \infty$, and
$$
  {L(r) \over Q(r)}
  = {k \over r\, \sqrt{\Lambda \gamma} }
  - {\gamma^{\,\prime}(r) \over 4 \sqrt{\Lambda}\, \gamma(r)^{3/2}}
  \rightarrow 0 \qquad (r \rightarrow \infty);
$$
moreover, $L/Q \in BV(\Any,\infty)$ as a result of Lemma 1, since
$$
  \left({1 \over r\sqrt{\gamma}}\right)'
  = - {\gamma^{\,\prime} \over 2 r \gamma^{3/2}} - {1 \over r^2
\sqrt{\gamma}}
  \in L^1(\Any,\infty).
$$
 
In particular, by Proposition 4 in the Appendix, there is $r_0 > 0$ such
that
$2 \modul{L(r)} \le Q(r)$ $(r \ge r_0)$ and $L/(Q - L) \in BV([r_0,
\infty))$.
Therefore $L$ and $Q$ satisfy conditions (C1), (C2) and (C3)' of
Proposition
2 with $M \equiv 0$.

Corollary 3 shows that there is a constant $C > 0$ such that
$$
  {\modul{v(r_0)}^2 \over C} \le \modul{v(r)}^2 \le C\,\modul{v(r_0)}^2,
\qquad
  {\modul{w(r_0)}^2 \over C} \le \modul{w(r)}^2 \le C\,\modul{w(r_0)}^2,
$$
i.e. that $v$ and $w$ are of the same size;
we now use this estimate to prove that
$$
  \liminf_{r\rightarrow\infty} {\int_{r_0}^r (u_1^2 + u_2^2)
      \over \int_{r_0}^r (y_1^2 + y_2^2)} > 0.
$$
To this end, we study the oscillation behaviour of $v$ by means of
the Pr\"ufer transformation.

There are locally absolutely continuous functions
$\varrho : [r_0, \infty) \rightarrow (0, \infty)$, and
$\vartheta : [r_0, \infty) \rightarrow \R$, such that
$\vartheta(r_0) \in [-3\pi/4, 5\pi/4)$, and
$$
  v = \varrho\,\pmatrix{\cos\vartheta \cr \sin\vartheta \cr};
$$
$\vartheta$ is a solution of the differential equation
$$
  \vartheta'(r) = L(r) \sin 2 \vartheta(r) + Q(r).
$$
Introducing a new independent variable by the transformation
$s(r) := \int_{r_0}^r Q$,
and setting $\Theta(s(r)) = \vartheta(r)$,
we have
$$
  \Theta'(s(r)) = 1 + {L(r) \over Q(r)}\,\sin 2 \Theta(s(r)).
$$
From the definition of $r_0$, we find that
$\Theta'(s) \in [1/2, 3/2]$ $(s \ge 0)$.

For $n \in \N$ we define
$$\eqalign{
  J_n &:= \{ s \ge 0 \mathrel{|} \Theta(s) \in [-3\pi/4, -\pi/4] + n\pi \},
\cr
  K_n &:= \{ s \ge 0 \mathrel{|} \Theta(s) \in [-\pi/4, \pi/4] + n\pi \};
\cr}$$
then
${\displaystyle \qquad
  {\pi \over 3} \le \modul{J_n},\modul{K_n} \le \pi
  \qquad (n \in \N_0)
}$

($\modul{\Any}$ here denotes the length of the interval).
Adjusting $r_0$ if necessary, we can assume that $0 = \inf J_1$.

Observing that $\gamma > \Lambda$, and $Q = \sqrt{\gamma\Lambda}$,
we find for $r \ge r_0$
$$\eqalign{
  {\int_{r_0}^r (u_1^2 + u_2^2)
     \over \int_{r_0}^r (y_1^2 + y_2^2)}
  &=
  {\int_{r_0}^r (\sqrt{\Lambda/\gamma}\, v_1^2 + \sqrt{\gamma/\Lambda}\,
v_2^2)
     \over \int_{r_0}^r (\sqrt{\Lambda/\gamma}\, w_1^2 +
     \sqrt{\gamma/\Lambda}\, w_2^2)}
  \ge
  {\int_{r_0}^r \sqrt{\gamma/\Lambda}\, v_2^2
     \over \int_{r_0}^r \sqrt{\gamma/\Lambda}\,(w_1^2 + w_2^2)}
\cr
  &=
  {\int_0^s V_2^2 \over \int_0^s (W_1^2 + W_2^2)}
  \ge
  {\int_0^s V_2^2 \over C\,\modul{w(r_0)}^2\,s(r)},
\cr}$$
where $v(r) = V(s(r))$, $w(r) = W(s(r))$ and $C$ is the
constant from Corollary 3.

Now setting $N(r) := \max \{n\in\N \mathrel{|} K_n \subset [0, s(r)]\}$
and noting that on $J_n$,
$$
  V_2^2 \ge {V_1^2 + V_2^2 \over 2} \ge {\modul{v(r_0)}^2 \over 2 C},
$$
we conclude that
$$\eqalignno{
  {\int_{r_0}^r (u_1^2 + u_2^2)
     \over \int_{r_0}^r (y_1^2 + y_2^2)}
  &\ge
  {\sum_{n=1}^{N(r)} \int_{J_n} V_2^2 \over C\,\modul{w(r_0)}^2\,s(r)}
  \ge
  {\pi\,\modul{v(r_0)}^2\,N(r) \over 6 C^2\,\modul{w(r_0)}^2
     \sum_{n=1}^{N(r) + 1} (\modul{J_n} + \modul{K_n})} &
\cr
  &\ge
  {\modul{v(r_0)}^2 \over 12 C^2\,\modul{w(r_0)}^2} \cdot
  {N(r) \over N(r) + 1}
  \rightarrow
  {\modul{v(r_0)}^2 \over 12 C^2\,\modul{w(r_0)}^2} > 0
    \qquad (r\rightarrow \infty). & \square
\cr}$$

\bigskip
{\bf Appendix.}
\smallskip
The great advantage of regularity assumptions on coefficient functions
in terms of their differentiability is based on the existence of a linear
and multiplicative differential calculus, by which subsequent estimates
of terms containing these coefficient functions are conveniently
accessible.
Unfortunately, such a calculus does not exist for the variation of
functions
on the real line, which generally behaves like the integral of the absolute
value of the derivative of the function, but cannot be treated as such
unless the function is differentiable or at least locally absolutely
continuous.
The approximation of functions of bounded variation by differentiable
functions, as in Weidmann [16] p. 368, is of limited scope
(and can hardly be effected in conditions like our (B2)),
whereas going back to the definition of the variation,
$$
  \mathop{\rm Var}_I f := \sup \sum_{j=1}^N \modul{f(x_j) - f(x_{j-1})}
$$
(where the supremum ranges over all finite partitions
$a = x_0 < x_1 < \dots < x_N = b$, $N\in\N$,
of the interval $I = [a, b]$),
is tedious and tends to obscure the line of reasoning.

On the other hand, it appears unsatisfactory to assume differentiability
of coefficients in situations where only their variation, but not their
derivatives occur naturally.

Therefore we have collected in this appendix some properties of functions
of bounded variation which have served as a substitute for the differential
calculus in the main body of this paper,
so that any approximation or reference to the definition above could be
avoided there, while working with minimal regularity requirements.
In particular, Proposition 3 is interesting as a weaker surrogate for the
product rule, and is frequently used in both Sections 2 and 3.
Similarly, Proposition 4 replaces the quotient rule in several instances
in our proofs.
The outlook of Proposition 5 is more restricted to the purpose of our
paper; it shows that condition (A3) of Theorem 1 is valid for all
$\lambda \in \R$ if it holds for two distinct values of $\lambda$.

\medskip
{\bf Proposition 3.}\qquad
{\it
Let $I \subset \R$ be an interval and $g \in BV(I)$,
$f\in AC\loc(I) \cap L^\infty(I)$.
Then
$$
  \mathop{\rm Var}_I f g
     \le \int_I \modul{f' g} + \norm{f}_\infty\mathop{\rm Var}_I g;
$$
in particular,
$f' g \in L^1(I)$ implies $f g \in BV(I)$.
}

\medskip
{\it Proof.}\qquad
Let $x, y \in I$, $x < y$.
By the Jordan decomposition theorem (Fichtenholz [5] \S 570),
there are non-decreasing functions $g_+, g_- : I \rightarrow \R$
such that $g = g_+ - g_-$, and
$\mathop{\rm Var}_{[x,y]} g = g_+(y) + g_-(y) - g_+(x) - g_-(x)$.

The formula for integration by parts for Stieltjes integrals
gives
$$
  f(y) g(y) - f(x) g(x)
  = \int_x^y f' g + \int_x^y f\,dg_+ - \int_x^y f\,dg_-.
$$
The mean value theorem for Stieltjes integrals (Fichtenholz [5] \S 582)
implies that there exist $z_+, z_- \in [x, y]$ such that
$$
  \int_x^y f\,dg_\pm = f(z_\pm)\,(g_\pm(y) - g_\pm(x));
$$
thus
$$
  \modul{f(y) g(y) - f(x) g(x)}
  \le \int_x^y \modul{f' g} + \norm{f}_\infty\,(g_+(y) - g_+(x))
    + \norm{f}_\infty\,(g_-(y) - g_-(x)).
$$
Applying this estimate to the individual intervals of partitions of $I$,
we obtain the assertion.
\hfill $\square$

\medskip
{\it Remark 7.}\qquad
The simpler estimate
$\mathop{\rm Var}_I f g \le
\norm{f}_\infty \mathop{\rm Var}_I g + \norm{g}_\infty \mathop{\rm Var}_I
f$
(for $f, g \in BV(I)$) is tempting;
but using this instead of Proposition 3 in the final step of the proof of
Theorem 1 requires the slightly stronger condition (A4)' (cf. Remark 6)
instead of (A4).

\medskip
{\bf Proposition 4.}\qquad
{\it
Let $I \subset \R$ be an interval and $f, g : I \rightarrow \R$.
Assume $g > 0$ and $\varepsilon := \sup_I \modul{f} / g < 1$.
Then
$$
  (1 - \varepsilon)^2\,\mathop{\rm Var}_I {f \over g - f}
  \le
  \mathop{\rm Var}_I {f \over g}
  \le
  (1 + \varepsilon)^2\,\mathop{\rm Var}_I {f \over g - f};
$$
in particular, $f/g \in BV(I)$ if and only if $f/(g-f) \in BV(I)$.
}

\medskip
{\it Proof.}\qquad
For $x, y \in I$ we have
$$
  \Modul{{f(y) \over g(y) - f(y)} - {f(x) \over g(x) - f(x)}}
  = {g(x) \over g(x) - f(x)} \cdot {g(y) \over g(y) - f(y)} \cdot
  \Modul{{f(y) \over g(y)} - {f(x) \over g(x)}},
$$
and the assertion follows observing that
$$
  {1 \over 1 + \varepsilon} \le {g \over g - f} \le {1 \over 1 -
\varepsilon}.
\eqno \square
$$

\medskip
{\bf Proposition 5.}\qquad
{\it
Let $a > 0$, $m, q : (a,\infty) \rightarrow \R$ functions such that
$$
  \lim_{r\rightarrow\infty} q(r) = \infty,
\qquad
  \limsup_{r\rightarrow\infty} {\modul{m(r)} \over q(r)} < 1.
$$
Then $\{\lambda \in \R \mathrel{|} {\displaystyle {m\over q - \lambda}} \in
BV(\Any,\infty)\}$
is either empty, or the whole real line, or has exactly one element.
}

\medskip
{\it Proof.}
Assume this set has two distinct elements $\lambda_1, \lambda_2$; we then
show that it is all of $\R$.
Let $\lambda\in\R$.
Then there is $r_0 > 0$ such that $m/(q - \lambda_j) \in BV[r_0,\infty)$,
and
$$
  q(r) - \lambda \ge {q(r) - \lambda_j \over 2} \ge 1
  \qquad (r \ge r_0, j \in \{1, 2\}).
$$

Defining $\mu_j := {\displaystyle {\lambda_j - \lambda \over \lambda_2 -
\lambda_1}}$
$(j \in \{1,2\})$, and noting that
$1 = \mu_2 - \mu_1$, $\lambda = \lambda_1\,\mu_2 - \lambda_2\,\mu_1$,
we find
$$\eqalign{
  &\Modul{{m(y) \over q(y) - \lambda} - {m(x) \over q(x) - \lambda}}
  = {\modul{m(y) q(x) - \lambda m(y) - m(x) q(x) - \lambda m(x)} \over
    (q(y) - \lambda)\,(q(x) - \lambda)}
\cr
  &\qquad\qquad\qquad
  \le 4 \modul{\mu_2}\,
    \Modul{{m(y) \over q(y) - \lambda_1} - {m(x) \over q(x) - \lambda_1}}
  + 4 \modul{\mu_1}\,
    \Modul{{m(y) \over q(y) - \lambda_2} - {m(x) \over q(x) - \lambda_2}}
\cr}$$
for $x, y \in [r_0, \infty)$.
Thus we obtain
$$
  \mathop{\rm Var}_{[r_0, \infty)} {m \over q - \lambda}
  \le 4 \modul{\mu_2}\,\mathop{\rm Var}_{[r_0,\infty)} {m \over q -
\lambda_1}
  + 4 \modul{\mu_1}\,\mathop{\rm Var}_{[r_0,\infty)} {m \over q -
\lambda_2}
  < \infty.
\eqno \square
$$

\medskip
{\it Remark 8.}\qquad
If $m$ is constant, the $\lambda\,m$ terms in the numerator cancel, and
it is easy to see that $m/(q-\lambda)$ is of bounded variation either for
every, or for no real value of $\lambda$ (as observed in [12]).

\medskip
\baselineskip=12pt
{\bf Acknowledgement.}\quad
The greater part of this work has been done during the 
visit of the second author at Munich University. He would like to express
his gratitude to Professors H. Kalf, E. Wienholtz and A. M. Hinz for their
hospitality.
  
\medskip
{\bf References.}

\leftskip=.4in
\parindent=-.4in

[1] Behncke H, Absolute continuity of Hamiltonians with von Neumann Wigner
\break potentials. {\it Proc. Amer. Math. Soc.,} {\bf 111} (1991),
373--384.

[2] Dunford N, Schwartz JT, {\it Linear operators II.} Wiley, New York,
1963.

[3] Erd\'elyi A, Note on a paper by Titchmarsh. {\it Quart. J. Math. Oxford
(2),} {\bf 14} (1963), 147--152.

[4] Evans WD, Harris BJ, Bounds for the point spectra of separated Dirac
operators. {\it Proc. Roy. Soc. Edinburgh,} {\bf 88A} (1981), 1--15.

[5] Fichtenholz GM, {\it Differential- und Integralrechnung III.} Deutscher
Verlag der \break Wissenschaften, Berlin, 1973.

[6] Gilbert DJ, Pearson DB, On subordinacy and analysis of the spectrum of
one-dimensional Schr\"odinger operators. {\it J. Math. Anal. Appl.,} {\bf
128}
(1987), 30--56.

[7] Gilbert DJ, On subordinacy and analysis of the spectrum of
Schr\"odinger
operators with two singular endpoints. {\it Proc. Roy. Soc. Edinburgh,}
{\bf 112A} (1989), 213--229.

[8] Hinton DB, Shaw JK, Dirac systems with discrete spectra. {\it Can. J.
Math.,} {\bf 39} (1987), 100--122.

[9] Roos BW, Sangren WC, Spectra for a pair of singular first order
differential
equations. {\it Proc. Amer. Math. Soc.,} {\bf 12} (1961), 468--476.

[10] Rose ME, Newton RR, Properties of Dirac wave functions in a central
field.
{\it Phys. Rev. (2),} {\bf 82} (1951), 470--477.

[11] Schmidt KM, Dense point spectrum and absolutely continuous spectrum in
\break spherically symmetric Dirac operators. {\it Forum Math.,} {\bf 7}
(1995),
459--475.

[12] Schmidt KM, Absolutely continuous spectrum of Dirac systems with
potentials
infinite at infinity. {\it Math. Proc. Camb. Phil. Soc.,} {\bf 122} (1997),
377--384.

[13] Simon B, Bounded eigenfunctions and absolutely continuous spectra for
one-\break dimensional Schr\"odinger operators. {\it Proc. Amer. Math.
Soc.,} {\bf 124} (1996), 3361--3369.

[14] Thaller B, {\it The Dirac equation.} Springer, Berlin, 1992.

[15] Titchmarsh EC, On the nature of the spectrum in problems of
relativistic
quantum mechanics. {\it Quart. J. Math. Oxford (2),} {\bf 12} (1961),
227--240.

[16] Weidmann J, Oszillationsmethoden f\"ur Systeme gew\"ohnlicher
Differential- \break gleichungen. {\it Math. Z.,} {\bf 119} (1971),
349--373.

[17] Weidmann J, Absolutstetiges Spektrum bei Sturm-Liouville-Operatoren
und
\break  Diracsystemen. {\it Math. Z.,} {\bf 180} (1982), 423--427.

[18] Weidmann J, {\it Spectral theory of ordinary differential operators.}
Lect. Notes Math. {\bf 1258}, Springer, Berlin, 1987.

[19] Weidmann J, Uniform nonsubordinacy and the absolutely continuous
spectrum.
{\it Analysis}, {\bf 16} (1996), 89--99.

[20] Yamada O, On the spectrum of Dirac operators with the unbounded
potential
at infinity. {\it Hokkaido Math. J.,} {\bf 26} (1997), 439--449.

\leftskip=0in
\parindent=0in

\bye